\documentclass[12pt,reqno]{amsart}
\usepackage{amsmath,amsthm,hyperref, amssymb}

\newtheorem{theorem}{Theorem}

\newtheorem{proposition}[theorem]{Proposition}

\theoremstyle{remark}
\newtheorem*{remark}{Remark}
\theoremstyle{definition}

\numberwithin{theorem}{section} \numberwithin{equation}{section}
\numberwithin{example}{section}

\newcounter{eqnno}

\renewcommand{\mod}{\mathrm{mod}\hspace*{3pt}}

\numberwithin{equation}{section} 
 
 \address{Eureka, California, 95501.}
 \email{dean.hickerson@yahoo.com}
\address{School of Mathematics and Physics, University of Queensland, Brisbane 4072}
\address{Max-Planck-Institut f\"ur Mathematik, Vivatsgasse 7, 53111 Bonn.}
\email{etmortenson@gmail.com}

\title{Dyson's Ranks and Appell--Lerch Sums}

\dedicatory{\em{For Freeman Dyson in honour of his 90th birthday}}

\begin{document}

\setcounter{eqnno}{1}

\date{29 June 2014}

\author{Dean Hickerson and Eric Mortenson}



\begin{abstract}  
Denote by $p(n)$ the number of partitions of $n$ and by $N(a,M;n)$ the number of partitions of $n$ with rank congruent to $a$ modulo $M$.  We find and prove a general formula for Dyson's ranks by considering the deviation of the ranks from the average:
\begin{equation*}
D(a,M) := \sum_{n= 0}^{\infty}\left(N(a,M;n) - \frac{p(n)}{M}\right) q^n.
\end{equation*}
Using Appell--Lerch sum properties we decompose $D(a,M)$ into modular and mock modular parts so that the mock modular component is supported on certain arithmetic progressions, whose modulus we can control.  Using our decomposition, we show how our formula gives as a straightforward consequence Atkin and Swinnerton-Dyer's results on ranks as well as Bringmann, Ono, and Rhoades's results on Maass forms.  We also apply our techniques to a variation of Dyson's ranks due to Berkovitch and Garvan.

\end{abstract}
\maketitle

\setcounter{equation}{-1}

\setcounter{section}{-1}

\section{Definitions}
Let $q$ be a complex number with $0<|q|<1$ and define $\mathbb{C}^*:=\mathbb{C}-\{0\}$.  We have
{\allowdisplaybreaks \begin{gather}
(x)_n=(x;q)_n:=\prod_{i=0}^{n-1}(1-q^ix), \ \ (x)_{\infty}=(x;q)_{\infty}:=\prod_{i\ge 0}(1-q^ix),\notag\\
 j(x;q):=(x)_{\infty}(q/x)_{\infty}(q)_{\infty}=\sum_{n=-\infty}^{\infty}(-1)^nq^{\binom{n}{2}}x^n,\label{equation:theta-series}\\
 {\text{and }}\ \ j(x_1,x_2,\dots,x_n;q):=j(x_1;q)j(x_2;q)\cdots j(x_n;q).\notag
\end{gather}}%
where in the middle line the equivalence of product and sum follows from Jacobi's triple product identity.  Identity (\ref{equation:theta-series}) defines a theta function.  A {\em partial} theta function is half of a theta function, for example
\begin{equation}
\sum_{n=0}^{\infty}(-1)^nq^{\binom{n}{2}}x^n.\label{equation:partial-theta}
\end{equation}
Let $a$ and $m$ be integers with $m$ positive.  We further define
\begin{gather}
J_{a,m}:=j(q^a;q^m), \ \ \overline{J}_{a,m}:=j(-q^a;q^m), \ {\text{and }}J_m:=J_{m,3m}=\prod_{i\ge 1}(1-q^{mi}).
\end{gather}

%  A theta function with the wrong signs is a {\em false} theta function; 

\section{Introduction}\label{section:introduction}
Freeman Dyson conjectured a beautiful combinatorial description of Ramanujan's congruences for the partition function using a statistic which he called the rank.  The rank has further intrigue due to its connections with Ramanujan's mock theta functions.   Atkin and Swinnerton-Dyer proved Dyson's conjectures and gave additional examples of how the rank relates to theta and mock theta functions.  Later, Bringmann and Ono generalised Atkin and Swinnerton-Dyer's results by using the theory of harmonic weak Maass forms, a theory which generalises mock theta functions.  In this paper, we use Appell--Lerch sum properties to give an explicit formula for Dyson's ranks which then yields as straightforward consequence results of Atkin and Swinnerton-Dyer \cite{ASD} as well as Bringmann, Ono, Rhoades, and Zagier \cite{BrO2, BrOR, Za}.  

%``{\em I am extremely sorry for not writing you a single letter up to now$\dots$I discovered very interesting functions recently which I call `Mock' $\vartheta$-functions.  Unlike the `False' $\vartheta$-theta functions (studied partially by Prof. Rogers in his interesting paper) they enter mathematics as beautifully as the ordinary theta functions$\dots$}''

In Ramanujan's last letter to Hardy, he gave a list of seventeen functions which he called ``mock theta functions''  \cite{R1}.   Each mock theta function $f(q)$ was defined as a $q$-series, convergent  for $|q|<1$, such that for every root of unity $\zeta$, there is a theta function $\theta_{\zeta}(q)$ such that the difference $f(q)-\theta_{\zeta}(q)$ is bounded as $q\rightarrow \zeta$ radially; moreover, there is no single theta function which works for all $\zeta$ \cite{GOR}.   Ramanujan also stated identities relating mock theta functions to each other as well as to modular forms, which suggested that mock theta functions live inside of a vector space, which has a subspace consisting of modular forms.  Later, more mock theta identities were found in the Lost Notebook \cite{RLN, H1}.%, and among the new identities were the mock-theta conjectures \cite{H1}.

The Lost Notebook is the Rosetta Stone of $q$-series.  Indeed, numerous entries expand $q$-hypergeometric series in terms of theta functions (Rogers-Ramanujan type identities), Appell--Lerch sums (mock theta functions), or partial theta functions.   Partial theta functions  play roles in areas outside of number theory such as quantum invariants of $3$-manifolds \cite{LZ}.  Appell--Lerch sums appear in the context of black hole physics \cite{DMZ}.  Long-standing problems have been to determine the modularity of the mock thetas and to understand how various types of $q$-series representations relate to each other.

In Zwegers' breakthrough work \cite{Zw2}, he solved the modularity question for mock theta functions.    Although mock theta functions are not modular \cite{W3}, he showed that they can be completed to non-holomorphic functions which are modular.  As a result, mock theta functions may be viewed as holomorphic parts of harmonic weak Maass forms \cite{  BrO1, BrO2,BruF, Za}.

%In Duke et al \cite{Du}  we see that there are no Hecke eigenforms in the environment of harmonic Maass forms.  Rhoades \cite{R} demonstrated that Ramanujan's definition of mock theta function does not agree with  the theory of harmonic Maass forms definition \cite{Za}. 

\subsection{Dyson's ranks}

One way to study mock theta functions is through partitions.  A partition of a positive integer $n$ is a weakly-decreasing sequence of positive integers whose sum is $n$.    The partitions of $4$ are $(4),$ $(3,1),$ $(2,2),$ $(2,1,1),$ $(1,1,1,1)$.    We denote the number of partitions of $n$ by $p(n)$.  Among the most famous results in partitions are Ramanujan's congruences:
{\allowdisplaybreaks \begin{align*}
p(5n+4)\equiv &\  0 \pmod{5}, \\
p(7n+5)\equiv &\  0 \pmod{7} ,\\
p(11n+6)\equiv &\  0 \pmod{11},
\end{align*}}%
\noindent for which there are many proofs and generalisations.  To study partition congruences, one often constructs a function, called a statistic, which assigns an integer value to a partition.  Dyson \cite{D1} gave insight into the first two congruences with such a statistic, which he called the rank.  He defined the rank of a partition to be the largest part minus the number of parts.   In particular, the ranks of the five partitions of $4$ are $3,1,0,-1,-3,$ respectively, giving an equinumerous distribution of the partitions of $4$ into the five distinct residue classes mod $5$.  In general, one defines
\begin{equation*}
N(a,M;n) := \mathrm{number\ of\ partitions\ of\ } n \mathrm{\ with\ rank\ } \equiv a \ (\mathrm{mod\ } M),
\end{equation*}
which has the symmetry property $N(a,M;n)=N(M-a,M;n)$.  Dyson conjectured \cite{D1}, \cite[$(2.2)$--$(2.11)$]{ASD}:
\begin{subequations}
\begin{equation}
N(1,5;5n+1)=N(2,5;5n+1)\label{equation:ASD-2.2},
\end{equation}
\begin{equation}
N(0,5;5n+2)=N(2,5;5n+2),
\end{equation}
\begin{equation}
N(0,5;5n+4)=N(1,5;5n+4)=N(2,5;5n+4),\label{equation:Dyson-M5}
\end{equation}
\begin{equation}
N(2,7;7n)=N(3,7;7n),\label{equation:ASD-2.5}
\end{equation}
\begin{equation}
N(1,7;7n+1)=N(2,7;7n+1)=N(3,7,7n+1),\label{equation:ASD-2.6}
\end{equation}
\begin{equation}
N(0,7,7n+2)=N(3,7;7n+2),
\end{equation}
\begin{equation}
N(0,7;7n+3)=N(2,7;7n+3), \ \ N(1,7;7n+3)=N(3,7;7n+3),
\end{equation}
\begin{equation}
N(0,7;7n+4)=N(1,7;7n+4)=N(3,7;7n+4),
\end{equation}
\begin{equation}
N(0,7;7n+5)=N(1,7;7n+5)=N(2,7;7n+5)=N(3,7;7n+5),\label{equation:Dyson-M7}
\end{equation}
\begin{equation}
N(0,7;7n+6)+N(1,7;7n+6)=N(2,7;7n+6)+N(3,7;7n+6).\label{equation:ASD-2.11}
\end{equation}
\end{subequations}
In particular, indentities (\ref{equation:Dyson-M5}) and (\ref{equation:Dyson-M7}) with the symmetry property give the first two of Ramanujan's congruences.

%\footnote{ Atkin and Swinnerton-Dyer leave undefined the rank of the empty partition of $0$, but we define it to have rank $0$.}  

We define the general rank-difference
\begin{equation}
R(a,b,M,c,m):=\sum_{n= 0}^{\infty} \Big ( N(a,M;mn+c)-N(b,M;mn+c)\Big )q^n,\label{equation:R-def}
\end{equation}
where $a,b,c,m,M$ are integers with $0\le a,b<M$ and $0\le c<m$.  Dyson's conjectures (\ref{equation:ASD-2.2})--(\ref{equation:Dyson-M7}) can be written as rank-differences that are equal to zero.  Atkin and Swinnerton-Dyer not only proved Dyson's conjectures but also determined rank-differences \cite[Theorems $4$, $5$]{ASD} which are equal to theta and mock theta functions.  As an example of Atkin and Swinnerton-Dyer's other results, they proved for modulus $5$ \cite[Theorem $4$]{ASD}, slightly rewritten,
\begin{subequations}
\begin{equation}
R(1,2,5,0,5)=qg(q,q^5),\label{equation:ASD-6.11}
\end{equation}
\begin{equation}
R(0,2,5,0,5)+2R(1,2,5,0,5)={J_{1}J_{5}^3}/{J_{1,5}^3},\label{equation:ASD-6.12}
\end{equation}
\begin{equation}
R(0,2,5,1,5)={J_{5}^2}/{J_{1,5}},\label{equation:ASD-6.13}
\end{equation}
\begin{equation}
R(1,2,5,2,5)=J_{5}^2/J_{2,5},\label{equation:ASD-6.16}
\end{equation}
\begin{equation}
R(0,2,5,3,5)=-qg(q^2,q^5),\label{equation:ASD-6.17}
\end{equation}
\begin{equation}
R(0,1,5,3,5)+R(0,2,5,3,5)=J_{1}J_{5}^3/J_{2,5}^3,\label{equation:ASD-6.18}
\end{equation}
\end{subequations}
where
\begin{equation}
g(x,q):=x^{-1}\Big ( -1+\sum_{n= 0}^{\infty}\frac{q^{n^2}}{(x;q)_{n+1}(q/x;q)_n}\Big )\label{equation:g-def}
\end{equation}
is a universal mock theta function.

Many of Ramanujan's mock theta functions can be written in terms of rank-differences.  For two fifth order mock theta functions \cite{H1}:
{\allowdisplaybreaks\begin{subequations}
\begin{equation}
\chi_0(q):=\sum_{n= 0}^{\infty}\frac{q^n}{(q^{n+1})_n}=2+3qg(q,q^5)-\frac{J_5^2J_{2,5}}{J_{1,5}^2}=2+R(1,0,5,0,5),\label{equation:chi0}
\end{equation}
\begin{equation}
\chi_1(q):=\sum_{n= 0}^{\infty}\frac{q^n}{(q^{n+1})_{n+1}}=3qg(q^2,q^5)+\frac{J_5^2J_{1,5}}{J_{2,5}^2}=R(2,0,5,3,5)+R(2,1,5,3,5).\label{equation:chi1}
\end{equation}
\end{subequations}}%
%The seventh order functions may be written in terms of rank-differences \cite{H2}.  
Most results about rank-differences have $m=M$, see (\ref{equation:R-def}), but for the mock theta $f(q)$:
\begin{equation}
f(q):=\sum_{n= 0}^{\infty}\frac{q^{n^2}}{(-q)_n^2}=2-2g(-1,q)=R(0,1,2,0,1).\label{equation:3rd-f}
\end{equation}
So we see that the $n$-th Fourier coefficient of $f(q)$ is equal to the number of partitions of $n$ with even rank minus the number of partitions of $n$ with odd rank \cite{BrO1}.

\subsection{The Main Theorem:  a special case}  We give an example of how our Main Theorem \ref{theorem:Theorem} yields Dyson's conjectures and Atkin and Swinnerton-Dyer's results.  For integers $0\le a<M$, define
\begin{equation}
D(a,M) = D(a, M,q) := \sum_{n= 0}^{\infty}\left(N(a,M;n) - \frac{p(n)}{M}\right) q^n.\label{equation:D-def}
\end{equation}
In terms of our definition (\ref{equation:R-def})
\begin{equation}
R(a,b,M,0,1)=D(a,M)-D(b,M). \label{equation:R2D}
\end{equation}
The Main Theorem \ref{theorem:Theorem} decomposes (\ref{equation:D-def}) into modular and mock modular components.  For modulus $M=5$, the Main Theorem specialises to
{\allowdisplaybreaks \begin{align}
D(0,5)&=-2 q^5 g(q^5, q^{25})\label{equation:D05}\\
& \ \ \ \ + \frac{4}{5}\cdot \frac{J_5 J_{25}^3}{J_{5,25}^3}
+ \frac{4}{5}\cdot q^1\cdot \frac{J_{25}^2}{J_{5,25}}
- \frac{2}{5}\cdot q^2\cdot \frac{ J_{25}^2}{J_{10,25}}
+ \frac{2}{5}\cdot q^3\cdot \frac{J_{5} J_{25}^3}{J_{10,25}^3}\notag\\
D(1,5) &= D(4,5) = q^5 g(q^5, q^{25})  - q^8 g(q^{10}, q^{25})\label{equation:D15-D45}\\
& \ \ \ \  - \frac{1}{5} \cdot    \frac{J_{5} J_{25}^3}{J_{5,25}^3}
  - \frac{1}{5} \cdot q^1 \cdot \frac{  J_{25}^2}{J_{5,25}}
 + \frac{3}{5} \cdot q^2 \cdot \frac{ J_{25}^2}{J_{10,25}}
 - \frac{3}{5} \cdot q^3 \cdot \frac{ J_5 J_{25}^3}{J_{10,25}^3}\notag\\
D(2,5) &= D(3,5) = q^8 g(q^{10}, q^{25}) \label{equation:D25-D35}\\
& \ \ \ \  - \frac{1}{5}\cdot \frac{J_5 J_{25}^3}{J_{5,25}^3}
- \frac{1}{5}\cdot q^1 \cdot \frac{ J_{25}^2}{J_{5,25}}
- \frac{2}{5}\cdot q^2 \cdot \frac{ J_{25}^2}{J_{10,25}}
+ \frac{2}{5}\cdot q^3 \cdot \frac{ J_5 J_{25}^3}{J_{10,25}^3}.\notag
\end{align}}%
The general form of the Main Theorem has twelve cases depending on $a$ and $M$ modulo $3$.

Identities (\ref{equation:ASD-2.2})--(\ref{equation:Dyson-M5}) and (\ref{equation:ASD-6.11})--(\ref{equation:ASD-6.18}) are now immediate.  Because of definitions (\ref{equation:g-def}) and (\ref{equation:theta-series}), the terms in the Fourier expansions from the $g$'s and $J$'s will be in terms of powers of $q^5$.  Given that there are no $q^{5k+4}$ terms in the expansions  (\ref{equation:D05}) -- (\ref{equation:D25-D35}), it is easy to see the first of Ramanujan's congruences (\ref{equation:Dyson-M5}).  We also see that (\ref{equation:ASD-2.2}) follows from (\ref{equation:D15-D45}) and (\ref{equation:D25-D35}) since the  coefficients of the $g$ parts do not involve any powers of $q$ of the form $q^{5k+1}$ and the theta parts for those powers are both equal to $-1/5 \cdot q^1 \cdot J_{25}^2/J_{5,25}$.   We even see rank-differences (\ref{equation:chi0}) and (\ref{equation:chi1}).

To prove and use our Main Theorem, we will need a building block finer than the universal mock theta function $g(x,q)$.  Although all of the mock theta functions from Ramanujan's letter \cite{R1} can be expressed in terms of the universal mock theta function $g(x,q)$ \cite{W3, H1,H2}, many mock thetas found in the Lost Notebook \cite{RLN} cannot.  In \cite{HM}, we demonstrated that all of the classical mock theta functions \cite{R1, RLN} can be expressed in terms of Appell-Lerch sums \cite{L1} defined as follows.  Let $x$, $q$, and $z$ be nonzero complex numbers with $|q|<1$ and neither $z$ nor $zx$ equal to an integral power of $q$, then  
\begin{equation}
m(x,q,z):=\frac{1}{j(z;q)}\sum_{r=-\infty}^{\infty}\frac{(-1)^rq^{\binom{r}{2}}z^r}{1-q^{r-1}xz}.\label{equation:mxqz-def}
\end{equation}
From \cite[Theorem $2.2$]{H1} and \cite{HM}, we can express $g(x,q)$ in terms of $m(x,q,z):$
\begin{align}
g(x,q)&=-x^{-2}m(qx^{-3},q^3,x^3z)-x^{-1}m(q^2x^{-3},q^3,x^3z)+\frac{J_1^2j(xz;q)j(z;q^3)}{j(x;q)j(z;q)j(x^3z;q^3)},\label{equation:g-full}
\end{align}
where the right-hand side is actually $z$-independent.  %Appell--Lerch sums are finer building blocks than $g(x,q)$ in that they have special collapsing and expanding  properties \cite[Theorems $3.9$, $3.5$]{HM}.

The Main Theorem can be used to determine the arithmetic progressions $\mathcal{A}$ for which
\begin{equation}
\sum_{n\in \mathcal{A}, \ n\ge 0} \Big (N(a,M;n)-\frac{p(n)}{M}\Big )q^{n-{1}/{24}}\label{equation:Zag-def}
\end{equation}
is modular.  To sketch an example, take expansions (\ref{equation:D05}) -- (\ref{equation:D25-D35}).  The mock modular contributions are supported on arithmetic progressions of modulus $5$, e.g. for (\ref{equation:D05}) the mock contribution is supported on the $q$-terms of the form $q^{5n}$.  As a consequence, an arithmetic progression modulo $5$ from (\ref{equation:D05}) not containing the terms $q^{5n}$ will give a weakly holomorphic modular form.  Identity (\ref{equation:g-full}) and an Appell--Lerch sum property \cite[Theorems $3.5$]{HM}, which we recall in Section \ref{section:tech}, allows us to control the modulus of the progression on which the mock modular component is supported.  A detailed discussion is found in Section \ref{section:tech}.

\subsection{Harmonic weak Maass forms}
%\footnote{For her work on mock theta functions and Dyson's ranks \cite{BrO1, BrO2, BrOR}, Bringmann won the $2009$ Sastra Prize and the $2009$ Alfried Krupp F\"orderpreis.}

It turns out that our Main Theorem gives as a straightforward consequence recent celebrated results of Bringmann, Ono, and Rhoades \cite{BrO2, BrOR, Za}, in which they employ the theory of harmonic weak Maass forms in order to produce far-reaching generalizations of Dyson's and Atkin and Swinnerton-Dyer's observations on rank-differences. To give examples of their results we recall notation.  For an integer $t>0$, define $\mathfrak{f}_t:=\frac{2t}{\gcd (t,6)}$, $\mathfrak{l}_t:={\textup{lcm}}(2t^2,24)$, and $\tilde{\mathfrak{l}}_t:=\mathfrak{l}_t/24.$  Define the group $\Gamma_c$ by
\begin{equation}
\Gamma_c:=\Big < \Big ( \begin{matrix} 1 & 1 \\ 0 & 0 \end{matrix}\Big ),  \Big ( \begin{matrix} 1 & 0 \\ \mathfrak{l}_c^2 & 0 \end{matrix}\Big )\Big >.
\end{equation}
Bringmann and Ono showed that the deviation of the ranks from the average value is a mock theta function:

\begin{theorem}  {\cite[Theorem $1.3$]{BrO2}} If $0\le r <t$ are integers, then
\begin{equation}
\sum_{n=0}^{\infty}\Big (N(r,t;n)-\frac{p(n)}{t}\Big )q^{\mathfrak{l}_t n-\tfrac{\mathfrak{l}_t}{24}}\label{equation:BrO2-Thm1.3-eqn}
\end{equation}
is a holomorphic part of a weak Maass form of weight $1/2$ on  $\Gamma_t$.  Moreover, if $t$ is odd, then it is on $\Gamma_{1}\big ( 144 \ \mathfrak{f}_t^2 \ \tilde{\mathfrak{l}}_t\big ).$
\end{theorem}
\noindent They found arithmetic progressions of the deviation that are modular:
\begin{theorem}  {\cite[Theorem $1.4$]{BrO2}} If $0\le r <t$ are integers, where $t$ is odd, and $\mathcal{P}\nmid 6t$ is prime, then
\begin{equation}
\sum_{\substack{n\ge 1 \\ \big ( \frac{24\mathfrak{l}_t n -\mathfrak{l}_t}{\mathcal{P}}\big )=-\big( \frac{-24\tilde{\mathfrak{l}}_t}{\mathcal{P}}\big )}}\Big (N(r,t;n)-\frac{p(n)}{t}\Big )q^{\mathfrak{l}_t n-\tfrac{\mathfrak{l}_t}{24}}\label{equation:BrO2-Thm1.4-eqn}
\end{equation}
is a weight $1/2$ weakly holomorphic modular form on $\Gamma_{1}\big ( 144 \ \mathfrak{f}_t^2 \  \tilde{\mathfrak{l}}_t \ \mathcal{P}^4\big ).$
\end{theorem}
\begin{remark}
One could rewrite (\ref{equation:BrO2-Thm1.4-eqn}) as a sum over $n$ with $ \big ( \frac{1-24n}{\mathcal{P}}\big )=-1$.
\end{remark}

\noindent Bringmann, Ono,  and Rhoades, also found generalisations of Atkin and Swinnerton-Dyer's results on rank-differences, see also Zagier \cite{Za}:
\begin{theorem} {\cite[Theorem $1.1$]{BrOR}} Suppose that $t\ge 5$ is prime, $0\le r_1,r_2 < t$ and $0\le d < t$.  Then the following are true:
\begin{itemize}
\item[(1)] If $\big ( \frac{1-24d}{t} \big ) =-1$, then 
\begin{equation*}
\sum_{n=0}^{\infty} (N(r_1,t;tn+d)-N(r_2,t;tn+d))q^{24(tn+d)-1}
\end{equation*}
is a weight $1/2$ weakly holomorphic modular form on $\Gamma_1(576t^6).$
\item[(2)] Suppose that $\big ( \frac{1-24d}{t} \big ) =1$.  If $r_1,r_2\not \equiv \tfrac{1}{2}(\pm 1 \pm \alpha) \pmod t$, where $\alpha$ is any integer for which $0\le \alpha <2t$ and $1-24d\equiv \alpha ^2 \pmod {2t}$, then
\begin{equation*}
\sum_{n=0}^{\infty} (N(r_1,t;tn+d)-N(r_2,t;tn+d))q^{24(tn+d)-1}
\end{equation*}
is a weight $1/2$ weakly holomorphic modular form on $\Gamma_1(576t^6).$
\end{itemize}
\end{theorem}

%Also, \cite[Theorem $1.4$]{BrO2} yields as a consequence
%\begin{theorem} {\cite[Theorem $1.2$]{BrOR}} Suppose that $t>1$ is an odd integer.  If $0\le r_1, r_2 <t$ are integers, and $\mathcal{P} \not | \ 6t$ is prime, then
%\begin{equation*}
%\sum_{\substack{n\ge 1 \\ \big ( \frac{24\mathfrak{l}_t n -\mathfrak{l}_t}{\mathcal{P}}\big )=-\big( \frac{-24\tilde{\mathfrak{l}}_t}{\mathcal{P}}\big )}}(N(r_1,t;n)-N(r_2,t;n))q^{\mathfrak{l}_t n-\tfrac{\mathfrak{l}_t}{24}}
%\end{equation*}
%is a weight $1/2$ weakly holomorphic modular form on $\Gamma_{1}\big ( 144 \mathfrak{f}_t^2 \tilde{\mathfrak{l}}_t\mathcal{P}^4\big ).$
%\end{theorem}
%\noindent One sees that \cite[Theorems $1.2$]{BrOR} and \cite[Theorems $1.4$]{BrO2} sum over $n$ such that $\big ( \frac{1-24n}{\mathcal{P}}\big )=-1.$

To obtain their results,  Bringmann et al. first complete an expression similar to (\ref{equation:Zag-def}) to a weight $1/2$ harmonic weak Maass form.  They then use a generalisation of quadratic twists (see for example Proposition $22$ of \cite{Ko}) in order to obtain an expression in which the non-holomorphic terms are supported on a certain arithmetic progression.  Using quadratic twists again they eliminate the non-holomorphic terms to obtain that, say, (\ref{equation:BrO2-Thm1.4-eqn}) is a weakly holomorphic modular form.  Using orthogonality of Dirichlet characters they then refine their results to, for example, \cite[Theorem $1.1(1)$]{BrOR}.

\subsection{Results} 

We obtain more explicit results without appealing to the theory of harmonic weak Maass forms.  Where Bringmann et al first complete (\ref{equation:Zag-def}) to a weight $1/2$ harmonic weak Maass form, we use Appell--Lerch sum properties to decompose (\ref{equation:Zag-def}) into modular and mock modular components, similar to (\ref{equation:D05}) -- (\ref{equation:D25-D35}).   Thus the Main Theorem is stronger version of \cite[Theorem $1.3$]{BrO2}.  To construct an arithmetic progression of terms from (\ref{equation:Zag-def}) which is a weakly holomorphic modular form, one applies quadratic twists and orthogonality of Dirichlet characters to the modular component.  One then checks that the mock modular component of the decomposition is supported on a disjoint progression.

%Using classical methods and our results of \cite{HM}, we prove the Main Theorem \ref{theorem:Theorem}.  The Main Theorem gives an explicit form of {\cite[Theorem $1.3$]{BrO2}} in that we decompose (\ref{equation:D-def}) into mock modular and modular parts.  

In Section \ref{section:tech} we discuss the role Appell--Lerch sums and their properties \cite[Theorems $3.9$, $3.5$]{HM} play in the proof the Main Theorem and in its applications.  In Section \ref{section:theorem} we state the Main Theorem and give additional specializations relevant to Dyson's rank-differences \cite{ASD, D1}.  In Section \ref{section:proof} we prove the Main Theorem using classical methods.  In Section \ref{section:BrO} we demonstrate how the Main Theorem yields results of Bringmann et al \cite{BrO2, BrOR, Za} on Dyson's ranks and Maass forms.   In Section \ref{section:M2} we apply our techniques to the $M_2$ rank of Berkovitch and Garvan \cite{BG}, which is based on work of Dyson \cite{D2,D3}.  

\section{The Role of Appell--Lerch sums}\label{section:tech}

We first recall well-known Appell--Lerch sum properties.
 \begin{proposition}  \cite{L1, HM, Zw2} If $x$, $q$, $z$, $z_0$, and $z_1$ satisfy the appropriate conditions of definition (\ref{equation:mxqz-def}), then
{\allowdisplaybreaks \begin{subequations}
\begin{equation}
m(x,q,z)=m(x,q,qz),\label{equation:m-fnq-z}
\end{equation}
\begin{equation}
m(x,q,z)=x^{-1}m(x^{-1},q,z^{-1}),\label{equation:m-fnq-flip}
\end{equation}
\begin{equation}
m(qx,q,z)=1-xm(x,q,z),\label{equation:m-fnq-x}
\end{equation}
\begin{equation}
m(x,q,z_1)-m(x,q,z_0)=\frac{z_0J_1^3j(z_1/z_0;q)j(xz_0z_1;q)}{j(z_0;q)j(z_1;q)j(xz_0;q)j(xz_1;q)}.\label{equation:m-change-z}
\end{equation}
\end{subequations}}
\end{proposition}
In our work \cite{HM}, we introduced a heuristic relating partial theta functions to Appell--Lerch sums.  The heuristic suggests that identities involving partial theta functions have analogous identities in terms of Appell--Lerchs sums; one just needs to add a theta function \cite[Section $3$]{HM}, \cite{Mo1}.  So for example, trivial identities such as 
{\allowdisplaybreaks \begin{gather}
\sum_{n=0}^{\infty}(-1)^nx^nq^{\binom{n}{2}}-\sum_{n=0}^{\infty}(-1)^n(-x)^nq^{\binom{n}{2}}
=-2x\sum_{n=0}^{\infty}(-1)^n(-x^2q^3)^{n}q^{4\binom{n}{2}},\label{equation:analog-1}\\
\sum_{n=0}^{\infty}(-1)^nx^nq^{\binom{n}{2}}=\sum_{n=0}^{\infty}(-1)^n(-x^2q)^{n}q^{4\binom{n}{2}}-x\sum_{n=0}^{\infty}(-1)^n(-x^2q^3)^{n}q^{4\binom{n}{2}},\label{equation:analog-2}
\end{gather}}%
where we sum over roots of unity or break up the summation index modulo some number, e.g. $2$, have non-trivial analogs in terms of Appell--Lerch sums.  The analog of (\ref{equation:analog-1}) generalises to

\begin{theorem} \cite[Theorem $3.9$] {HM} \label{theorem:rootsof1}
Let $n$ and $k$ be integers with $0 \leq k < n$. Let $\omega$ be a primitive $n$-th root of unity. Then
\begin{align*}
\sum_{t=0}^{n-1}& \omega^{-kt} m(\omega^t x,q,z) =  n q^{-\binom{k+1}{2}} (-x)^k m\big(-q^{{\binom{n}{2}-nk}} (-x)^n, q^{n^2}, z' \big)+n\Psi_k^n(x,z,z^{\prime};q),%\label{equation:rootsof1}
\end{align*}
where
\begin{align}
\Psi_k^n&(x,z,z^{\prime};q):=\label{equation:msplit-theta}\\
&-\frac{ x^k z^{k+1} J_{n^2}^3} {j(z;q)j(z';q^{n^2})}\cdot\sum_{t=0}^{n-1} \frac{q^{{\binom{t+1}{2}+kt}} (-z)^t j\big(-q^{{\binom{n+1}{2}+nk+nt}} (-z)^n/z';q^{n^2}\big)
j(q^{nt} x^n z^n z';q^{n^2})}
{j\big(-q^{{\binom{n}{2}-nk}} (-x)^n z', q^{nt} x^n z^n;q^{n^2}\big)}.\notag
\end{align}
\end{theorem}

Let $\zeta_M:=e^{2\pi i /M}$, we have from \cite[$(3.13)$]{BrO2}:
\begin{align}
D(a,M)=\frac{1}{M}\sum_{j=1}^{M-1}\zeta_M^{-aj}\Big (1-\zeta_M^j\Big )\Big (1+\zeta_M^jg(\zeta_M^j,q)\Big ),\label{equation:D-to-g}
\end{align}
where the $g(x,q)$ is the universal mock theta function (\ref{equation:g-def}).  To decompose $D(a,M)$ into modular and mock modular components, we first expand (\ref{equation:D-to-g}) in terms of the $m(x,q,z)$ function, see (\ref{equation:g-full}).  Once expanded, we can use  Thereom \ref{theorem:rootsof1} to collapse the sum over roots of unity.  We then use (\ref{equation:g-full}) again to write the resulting $m(x,q,z)$ functions  in terms of the universal mock theta function $g(x,q)$.   This method gives our Main Theorem.

Producing generalisations of rank-differences such as  \cite[Theorem $1.4$]{BrO2} and \cite[Theorem $1.1$]{BrOR} requires another Appell--Lerch sum property from \cite{HM}.  Let us consider (\ref{equation:D25-D35}).  We remind the reader of the various ways of producing new modular forms from old ones such as by taking the twist by a Dirichlet character or by restricting the index of summation of the Fourier expansion to terms that lie in certain arithmetic progressions.  With this in mind we see that up to multiplication by a factor $q^{-1/24}$ that we can restrict the right-hand side (\ref{equation:D25-D35}) to the sequences $q^{5n}$, $q^{5n+1}$, or $q^{5n+2}$, and obtain a weakly holomorphic modular form.  In particular, we point out that these are the sequences to which $q^{8}g(q^{10},q^{25})$ does not contribute.  If we wish to consider sequences modulo $M$ instead of modulo $5$,  we rewrite the $q^{8}g(q^{10},q^{25})$ in terms of the $m(x,q,z)$ function, see (\ref{equation:g-full}), and use the following theorem with $n=M$ and $z^{\prime}=-1$:
\begin{theorem} \cite[Theorem $3.5$]{HM} \label{theorem:msplit-general-n} For generic $x,z,z'\in \mathbb{C}^*$ 
{\allowdisplaybreaks \begin{align*}
m(&x,q,z) = \sum_{r=0}^{n-1} q^{{-\binom{r+1}{2}}} (-x)^r m\big(-q^{{\binom{n}{2}-nr}} (-x)^n, q^{n^2}, z' \big)\\\notag
&+ \frac{z' J_n^3}{j(xz;q) j(z';q^{n^2})} \cdot \sum_{r=0}^{n-1}
\frac{q^{{\binom{r}{2}}} (-xz)^r
j\big(-q^{{\binom{n}{2}+r}} (-x)^n z z';q^n\big)
j(q^{nr} z^n/z';q^{n^2})}
{ j\big(-q^{{\binom{n}{2}}} (-x)^n z', q^r z;q^n\big )}.
\end{align*}}%
\end{theorem}
%\noindent Theorem \ref{theorem:msplit-general-n} is actually a consequence of Theorem \ref{theorem:rootsof1}.  
\noindent Instead of only allowing sequences module $5$ to which $q^8g(q^{10},q^{25})$ does not contribute, we now avoid contributions from terms of the form $q^{d_r}m(q^{Ma_r},q^{M^2b_r},-1)$.  Here $a_r$, $b_r$, and $d_r$ are integers.  By definition (\ref{equation:mxqz-def}) the Appell--Lech sums contribute only powers of $q^{M}$ to the Fourier expansion, so we just need to avoid $q$-exponents congruent to $d_r$ modulo $M$.

%We find and prove a general formula for $D(a,M)$, our main theorem,  and demonstrate how the main theorem is useful in understanding Dyson's ranks \cite{ASD,BrO2, BrOR,D1}.  

We list a few more facts which will be useful in the paper.  We note the $z=x^{-1}$ and $z=-x^{-3}$ specialisations of identity (\ref{equation:g-full}):
\begin{subequations}
\begin{equation}
g(x,q)=-x^{-1}m(q^2x^{-3},q^3,x^2)-x^{-2}m(qx^{-3},q^3,x^2),\label{equation:g-to-m}
\end{equation}
\begin{equation}
g(x,q)=-x^{-2}m(qx^{-3},q^3,-1)-x^{-1}m(q^2x^{-3},q^3,-1)+\Delta(x;q),\label{equation:g-spec1}
\end{equation}
\end{subequations}
where
\begin{align}
\Delta(x;q):=\frac{x^{-2}J_1J_{3}^3j(-x^{2};q)}{j(x;q)j(-qx^{3};q^3)j(-q^{2}x^{3};q^3)\overline{J}_{0,3}}.\label{equation:theta-delta}
\end{align}
For computing examples involving $g(x,q)$, we recall an identity {\cite[p. $32$]{RLN}, \cite[$(12.5.3)$]{ABI}}:
\begin{equation}
g(x,q)=-x^{-1}+qx^{-3}g(-qx^{-2},q^4)-qg(-qx^2,q^4)+\frac{J_2J_{2,4}^2}{xj(x;q)j(-qx^2;q^2)},\label{equation:g-split}
\end{equation}
as well as two properties:
\begin{align}
g(x,q)=g(q/x,q) {\text{ and }}  g(qx,q)=-x-x^2-x^3g(x,q).
\end{align}

\section{The Main Theorem}\label{section:theorem}
\begin{theorem}\label{theorem:Theorem}
If $0\leq a<M$, then
\begin{equation*}
D(a, M,q) = d(a, M,q) + T_{a, M}(q)
\end{equation*}
where
\newpage

\enlargethispage{3\baselineskip}
\begin{equation*}
d(a, M,q)=
 \begin{cases}
% M == 0 (mod 3)
%                             a == 0 (mod 3)
2\,m((-1)^{M+1}q^{M(M-1)/6},q^{M^2/3},z)
\hspace{.5in}\mathrm{if\ }a=0 \mathrm{\ and\ }M\equiv 0\ (\mod3)\\[.07in]
(-1)^a\ q^{-a(a+1)/6}\ m((-1)^{M+1}\ q^{M(M-2a-1)/6},q^{M^2/3},z)\\
\hspace{1in}+ (-1)^{a+1}\ q^{-a(a-1)/6}\ m((-1)^{M+1}\ q^{M(M-2a+1)/6},q^{M^2/3},z)\\
\hspace{1.5in}\mathrm{if\ }a\equiv 0\ (\mod3) \mathrm{,\ }M\equiv 0\ (\mod3)\mathrm{,\ and\ }a\neq0\\[.07in]
%                             a == 1 (mod 3)
(-1)^a\ q^{-a(a-1)/6}\ m((-1)^{M+1}\ q^{M(M-2a+1)/6},q^{M^2/3},z)\\
\hspace{1.5in}\mathrm{if\ }a\equiv 1\ (\mod3)\mathrm{\ and\ }M\equiv 0\ (\mod3)\\[.07in]
%                             a == 2 (mod 3)
%%(-1)^{a+1}\ q^{-a(a+1)/6}\ m((-1)^{M+1}\ q^{M(M-2a-1)/6},q^{M^2/3},z)\\
%%\hspace{1.5in}\mathrm{if\ }a\equiv 2\ (\mod3)\mathrm{\ and\ }M\equiv 0\ (\mod3)\\
\rule{5in}{.5pt}\\
% M == 1 (mod 3)
%                             a == 0 (mod 3)
0\hspace{1.44in}\mathrm{if\ }a=0 \mathrm{\ and\ }M=1\\[.07in]
2+2\ (-1)^{M+1}\ q^{M(M-1)/6}\ g((-1)^{M+1}\ q^{M(M-1)/6},q^{M^2})\\
\hspace{1.5in}\mathrm{if\ }a=0 \mathrm{,\ }M\equiv 1\ (\mod3)\mathrm{,\ and\ }M\neq1\\[.07in]
(-1)^a q^{-a(a+1)/6}
 +(-1)^{M+a+1}\ q^{(M(M-2a-1)-a(a+1))/6}\ g((-1)^{M+1} q^{M(M-2a-1)/6},q^{M^2})\\
\hspace{1in}+\;(-1)^{M+a+1}\ q^{(M(M+2a-1)-a(a-1))/6}\ g((-1)^{M+1} q^{M(M+2a-1)/6},q^{M^2})\\
\hspace{1.5in}\mathrm{if\ }a\equiv 0\ (\mod3)\mathrm{,\ }M\equiv 1\ (\mod3) \mathrm{,\ }a\neq0 \mathrm{,\ and\ }M\neq 2a+1\\[.07in]
(-1)^a q^{a(7a+5)/6}\ g(q^{2aM/3},q^{M^2})
\hspace{.5in}\mathrm{if\ }a\equiv 0\ (\mod3)\mathrm{,\ }a\neq0 \mathrm{,\ and\ }M=2a+1\\[.07in]
%                             a == 1 (mod 3)
%(-1)^a\ q^{-a(a-1)/6}
%+(-1)^{M+a+1}\ q^{(M(M-2a+1)-a(a-1))/6}\ g((-1)^{M+1}\ q^{M(M-2a+1)/6},q^{M^2})\\
%\hspace{1.05in}+\;(-1)^{a+1}\ q^{(2M^2-a(a+1))/6}\ g((-1)^{M+1}\ q^{M(3M-2a-1)/6},q^{M^2})\\
%\hspace{1.5in}\mathrm{if\ }a\equiv 1\ (\mod3)\mathrm{,\ }M\equiv 1\ (\mod3)\mathrm{,\ and\ }M\neq 2a-1\\[.07in]

%(-1)^{a+1}\ q^{(a-1)(7a-2)/6}\ g(q^{2(a-1)M/3},q^{M^2})
%\hspace{.5in}\mathrm{if\ }a\equiv 1\ (\mod3)\mathrm{\ and\ }M=2a-1\\[.07in]

%                             a == 2 (mod 3)
(-1)^{M+a+1}\ q^{(M(M+2a+1)-a(a+1))/6}\ g((-1)^{M+1}\ q^{M(M+2a+1)/6},q^{M^2})\\
\hspace{1in}+ (-1)^{a+1}\ q^{(2M^2-a(a-1))/6}\ g((-1)^{M+1}\ q^{M(3M+2a-1)/6},q^{M^2})\\
\hspace{1.5in}\mathrm{if\ }a\equiv 2\ (\mod3)\mathrm{\ and\ }M\equiv 1\ (\mod3)\\
\rule{5in}{.5pt}\\
% M == 2 (mod 3)
%                             a == 0 (mod 3)
2\ (-1)^M\ q^{M(M+1)/6}\ g((-1)^{M+1}\ q^{M(M+1)/6},q^{M^2})\\
\hspace{1.5in}\mathrm{if\ }a=0 \mathrm{\ and\ }M\equiv 2\ (\mod3)\\[.07in]
(-1)^{a+1}\ q^{-a(a-1)/6}
 + (-1)^{M+a}\ q^{(M(M-2a+1)-a(a-1))/6}\ g((-1)^{M+1}\ q^{M(M-2a+1)/6},q^{M^2})\\
\hspace{1.2in}+\;(-1)^{M+a}\ q^{(M(M+2a+1)-a(a+1))/6}\ g((-1)^{M+1}\ q^{M(M+2a+1)/6},q^{M^2})\\
\hspace{1.5in}\mathrm{if\ }a\equiv 0\ (\mod3) \mathrm{,\ }M\equiv 2\ (\mod3) \mathrm{,\ }a\neq0\mathrm{,\ and\ }M\neq 2a-1\\[.07in]
(-1)^{a+1}\ q^{a(7a-5)/6}\ g(q^{2aM/3},q^{M^2})
\hspace{.5in}\mathrm{if\ }a\equiv 0\ (\mod3)\mathrm{\ and\ }M=2a-1\\
%                             a == 1 (mod 3)
(-1)^{M+a}\ q^{(M(M+2a-1)-a(a-1))/6}\ g((-1)^{M+1}\ q^{M(M+2a-1)/6},q^{M^2})\\
\hspace{1in}+ (-1)^a\ q^{(2M^2-a(a+1))/6}\ g((-1)^{M+1}\ q^{M(3M-2a-1)/6},q^{M^2})\\
\hspace{1.5in}\mathrm{if\ }a\equiv 1\ (\mod3)\mathrm{\ and\ }M\equiv 2\ (\mod3)
%                             a == 2 (mod 3)
%%(-1)^{a+1}\ q^{-a(a+1)/6}
%%+ (-1)^{M+a}\ q^{(M(M-2a-1)-a(a+1))/6}\ g((-1)^{M+1}\ q^{M(M-2a-1)/6},q^{M^2})\\
%%\hspace{1.19in}+ (-1)^a\ q^{(2M^2-a(a-1))/6}\ g((-1)^{M+1}\ q^{M(3M+2a-1)/6},q^{M^2})\\
%%\hspace{1.5in}\mathrm{if\ }a\equiv 2\ (\mod3)\mathrm{,\ }M\equiv 2\ (\mod3)\mathrm{,\ and\ }M\neq 2a+1\\
%%(-1)^a\ q^{(a+1)(7a+2)/6}\ g(q^{2(a+1)M/3},q^{M^2})
%%\hspace{.5in}\mathrm{if\ }a\equiv 2\ (\mod3)\mathrm{\ and\ }M=2a+1\\
\end{cases}
\end{equation*}
and $T_{a,M}(q)$ is a theta function.   For $M\equiv 0 \pmod 3$, $T_{a,M}(q)$ also depends on the $z$.
\end{theorem}

We have removed the cases determined by $D(a,M)=D(M-a,M).$  We note that the case $a=0$ is almost a special case of $a\equiv0 \pmod3$; we just have to add $1$.  This is related to the empty partition of $0$, whose rank Atkin and Swinnerton-Dyer leave undefined.  We have defined it to have rank $0$.

For a theta function $f(q)$ (and for many other functions) there is an associated fractional exponent $b/a$ such that $q^{b/a}\cdot f(q)$ is, in some ways, simpler than $f(q).$  Modular properties of theta functions are easier to state when you multiply by $q^{b/a}.$  Let us call $b/a$ the multiplier exponent or $\lambda$.  The $\lambda$ for both $J_{a,m}$ and $\overline{J}_{a,m}$ is $(m-2a)^2/8m.$  For $m(q^a,q^m,z)$, $\lambda$ is $-(m-2a)^2/8m.$  We note that the $\lambda$'s for all of the cases and summands therein of the theorem are equal to $-1/24.$  When multiplied by $q^{-1/24}$,  the term $T_{a,M}(q)$ becomes a weight $1/2$ weakly holomorphic modular form.  

 In Section \ref{section:proof},  we see that we can write $T_{a,M}(q)$ in terms of the theta function $j(x;q)$.   We explain the arbitrary $z$ as well as how and when $T_{a,M}(q)$ depends on it.  When using Theorem \ref{theorem:rootsof1}, we have an arbitrary $z^{\prime}$ in both the resulting $m(x,q,z)$ function and the quotients of theta functions.  For $M\equiv 1,  \ 2 \pmod{3}$ we specialise the $z^{\prime}$ in order to form $g(x,q)$ terms as in (\ref{equation:g-full}).  For $M\equiv 0 \pmod 3$, we cannot obtain $g(x,q)$ terms, so we leave the $z^{\prime}$ arbitrary for the reader to decide; this situation is not unlike identity (\ref{equation:g-full}) in which the right-hand side is $z$-independent.

\section{Examples of Formulas for $D(a,M)$}
We state formulas for small values of $M$ and give examples of proofs for some of the individual $D(a,M,q)$'s.
\subsection{ Case {M=2}:}
\begin{align}
D(0,2)&=2qg(-q,q^4) + \frac{1}{2}\cdot \frac{J_{1,2}^2}{ J_{1}}\\
D(1,2)&=-D(0,2)=-2qg(-q,q^4) - \frac{1}{2}\cdot \frac{J_{1,2}^2}{  J_{1}}
\end{align}

To prove $D(0,2)$, we use (\ref{equation:g-split}) to obtain
\begin{align*}
g(-1,q)&=1-2qg(-q,q^4)-\frac{1}{2}\cdot \frac{J_1^3}{J_2^2}.
\end{align*}
Using (\ref{equation:D-to-g}) we arrive at
\begin{align*}
D(0,2)&=\frac{1}{2}\cdot 1\cdot(1-(-1))\cdot (1-g(-1,q))=2qg(-q,q^4)+\frac{1}{2}\cdot \frac{J_1^3}{J_2^2}.
\end{align*}

\subsection{Case M=3:}
\begin{align}
D(0,3)&=2m(q,q^3,-1) - \frac{1}{3}\cdot \frac{J_{1,2}^2 }{\overline{J}_{1,3}}\\
D(1,3)&=D(2,3) = -D(0,3)/2 = -m(q,q^3,-1) + \frac{1}{6}\cdot \frac{J_{1,2}^2}{\overline{J}_{1,3}}
\end{align}
%\subsection{Case M=4:}
%\begin{align}
%D(0,4)&= 2-2q^2g(-q^2,q^{16})+\frac{\overline{J}_{1,2}^2-4\overline{J}_{4,8}^2-2J_{1,2}J_{2,4}}{4J_1}\\
%D(1,4)=D(3,4)&= -1+q^2g(-q^2,q^{16})+q^5g(-q^6,q^{16})+\frac{4\overline{J}_{4,8}^2-\overline{J}_{1,2}^2}{4J_1}\\
%D(2,4)&=-2q^5g(-q^6,q^{16})+\frac{\overline{J}_{1,2}^2-4\overline{J}_{4,8}^2+2J_{1,2}J_{2,4}}{4J_1}
%\end{align}

For $D(0,3)$, we specialize $M=3$ and $z^{\prime}=-1$ in (\ref{equation:Dfull-a-equal-0-and-3-div-M}) and note straightforward facts such as $j(\zeta_3;q)=(1-\zeta_3)J_3$ and $j(-\zeta_3;q)=(1+\zeta_3)J_1^2J_6/J_2J_3$. It follows that
\begin{align*}
D(0,3)&=2m(q,q^3,-1)+\frac{1}{3}\frac{J_1^2}{\overline{J}_{0,1}}\Big [ \frac{(\zeta_3-\zeta_3^2)}{(1-\zeta_3)J_3}\frac{(1+\zeta_3)J_1^2J_6}{J_2J_3} 
+  \frac{(\zeta_3^2-\zeta_3)}{(1-\zeta_3^2)J_3}\frac{(1+\zeta_3^2)J_1^2J_6}{J_2J_3}  \Big ]\\
&=2m(q,q^3,-1)-\frac{J_1^5J_6}{3J_2^3J_3^2},
\end{align*}
where the last line follows from  the fact $1+\zeta_3+\zeta_3^2=0$ and some simplifying.

\subsection{Case M=7:}
{\allowdisplaybreaks \begin{align}
&D(0,7) = 2 + 2 q^7 g(q^7, q^{49})\\
&         - \frac{8}{7}\cdot \frac{     J_{21,49}^2 }{ J_7}
+ \frac{6}{7}\cdot  q^1 \cdot \frac{ J_{49}^2}{  J_{7,49}}
- \frac{2}{7}\cdot q^2 \cdot \frac{J_{14,49}^2}{ J_{7}}
+ \frac{4}{7}\cdot  q^3 \cdot \frac{J_{49}^2 }{ J_{14,49}}
+ \frac{2}{7}\cdot q^4 \cdot \frac{J_{49}^2 }{ J_{21,49}}
- \frac{4}{7}\cdot q^6\cdot \frac{ J_{7,49}^2 }{ J_{7}}\notag\\
&D(1,7)= D(6,7) = -1 - q^7 g(q^7, q^{49}) + q^{16} g(q^{21}, q^{49})\label{equation:D17}\\
&  + \frac{6}{7}\cdot \frac{J_{21,49}^2}{ J_{7}}
- \frac{1}{7}\cdot q^1 \cdot \frac{J_{49}^2}{  J_{7,49}}
+ \frac{5}{7}\cdot q^2 \cdot \frac{J_{14,49}^2}{ J_{7}}
- \frac{3}{7}\cdot q^3 \cdot \frac{J_{49}^2 }{ J_{14,49}}
+ \frac{2}{7}\cdot q^4 \cdot \frac{J_{49}^2}{ J_{21,49}}
+ \frac{3}{7}\cdot q^6 \cdot \frac{J_{7,49}^2}{ J_7}\notag\\
&D(2,7) = D(5,7) = q^{13} g(q^{14}, q^{49}) - q^{16} g(q^{21}, q^{49})\\
& - \frac{1}{7} \cdot \frac{J_{21,49}^2}{ J_7}
- \frac{1}{7}\cdot  q^1 \cdot \frac{  J_{49}^2}{  J_{7,49}}
- \frac{2}{7}\cdot  q^2 \cdot \frac{J_{14,49}^2}{ J_{7}}
+ \frac{4}{7} \cdot q^3 \cdot \frac{J_{49}^2}{  J_{14,49}}
- \frac{5}{7}\cdot q^4 \cdot \frac{J_{49}^2}{  J_{21,49}}
+ \frac{3}{7} \cdot q^6 \cdot \frac{J_{7,49}^2}{ J_7}\notag\\
&D(3,7) = D(4,7) = - q^{13} g(q^{14}, q^{49})\\
&   - \frac{1}{7}\cdot \frac{     J_{21, 49}^2}{ J_{7}}
- \frac{1}{7}\cdot  q^1 \cdot \frac{   J_{49}^2 }{ J_{7, 49}}
- \frac{2}{7}\cdot  q^2 \cdot \frac {J_{14, 49}^2}{  J_{7}}
- \frac{3}{7}\cdot  q^3 \cdot \frac{J_{49}^2}{ J_{14, 49}}
+ \frac{2}{7}\cdot  q^4 \cdot \frac{J_{49}^2}{ J_{21, 49}}
- \frac{4}{7}\cdot  q^6 \cdot \frac{J_{7, 49}^2}{ J_{7}}\notag
\end{align}}%

Dyson's (\ref{equation:ASD-2.5})--(\ref{equation:ASD-2.11}) and Atkin and Swinnerton-Dyer's \cite[Theorem $5$]{ASD} are immediate. Dyson's conjecture (\ref{equation:Dyson-M7}) is clear as there are no terms of the form $q^{7n+5}$.  Considering coefficients of respective $g(x,q)$ terms easily leads to identities such as (\ref{equation:ASD-2.6}).  In particular, the $g$ expressions in $D(1,7)$ only contribute to the $q^{7n}$ and $q^{7n+2}$ progressions, so for the $q^{7n+1}$ terms we only have
\begin{equation*}
- \frac{1}{7}\cdot  q^1 \cdot \frac{  J_{49}^2}{  J_{7,49}}.
\end{equation*}

%\subsection{Case M=13:}

%In $(2.11)$ and $(2.12)$ of Andrews's smallest parts paper \cite{A8}, he's looking at the coefficients of  $q^{13k+6}$  in certain linear combinations of:
%\begin{align*}
%D(0,13)&\sim 2+2q^{26}g(q^{26},q^{169}),\\
%D(1,13)=D(12,13)&\sim -1 - q^{26} g(q^{26},q^{169}) + q^{56} g(q^{78},q^{169}),\\
%D(2,13)=D(11,13)&\sim q^{38} g(q^{39},q^{169}) - q^{56} g(q^{91},q^{169}),\\
%D(3,13)=D(10,13)&\sim -q^{-2} - q^{11} g(q^{13},q^{169}) - q^{38} g(q^{39},q^{169}),\\
%D(4,13)=D(9,13)&\sim q^{-2} + q^{11} g(q^{13},q^{169}) - q^{53} g(q^{65},q^{169}),\\
%D(5,13)=D(8,13)&\sim -q^{47} g(q^{52},q^{169}) + q^{53} g(q^{104},q^{169}),\\
%D(6,13)=D(7,13)&\sim q^{47} g(q^{52},q^{169}),
%\end{align*}
%where `$\sim$' means mod theta.  The exponents in the powers of $q$ that are multiplied by $g(q^{13k}, q^{169})$ are $11$, $26$, $38$, $47$, $53$, and $56$.  None of those are congruent to $6$ (mod $13$), so the $g(x,q)$ parts of the formulas for the rank generating functions do not contribute anything to the coefficients.  It turns out that the theta parts give multiples of $13$ for those two linear combinations \cite{A8, O}.  

%%%%%%%%%%%
%%%%%%%%%%%
\section{Proof of the Main Theorem}\label{section:proof}

The proofs of the twelve cases are similar, so we give only a few as examples.  For the sake of brevity, we do not write out an explicit theta function for each $T_{a,M}(q)$.  However, we point out that our arguments are effective in that we can easily keep track of the various summands of quotients of theta functions that arise from using (\ref{equation:m-change-z}) and Theorem \ref{theorem:rootsof1}.  

We recall that $\zeta_M:=e^{2\pi i /M}$, so by (\ref{equation:D-to-g}):
{\allowdisplaybreaks \begin{align}
&\sum_{n\geq 0}\left(N(a,M,n) - \frac{p(n)}{M}\right) q^n =\frac{1}{M}\sum_{j=1}^{M-1}\zeta_M^{-aj}\Big (1-\zeta_M^j\Big )\Big (1+\zeta_M^jg(\zeta_M^j,q)\Big )\notag\\
&=\frac{1}{M}\sum_{j=1}^{M-1}\zeta_M^{-aj}\Big (1-\zeta_M^j\Big )\Big (1-m\big (q^2\zeta_M^{-3j},q^3,-1\big )-\zeta_M^{-j}m\big (q\zeta_M^{-3j},q^3,-1\big ) + \zeta_M^j\Delta(\zeta_M^j;q) \Big )\notag\\
&= \frac{1}{M}\sum_{j=0}^{M-1}\zeta_M^{-aj}\Big (m\big (q\zeta_M^{3j},q^3,-1\big )-\zeta_M^{-j}m\big (q\zeta_M^{-3j},q^3,-1\big )
 -\zeta_M^j m\big (q\zeta_M^{3j},q^3,-1\big )\notag \\
 &\ \ \ \ \ \ \ \ \ \ +m\big (q\zeta_M^{-3j},q^3,-1\big )\Big )  +\frac{1}{M}\sum_{j=1}^{M-1}\zeta_M^{-aj}\Big (1-\zeta_M^j\Big ) \zeta_M^j\Delta(\zeta_M^j;q),\label{equation:DAMQ-heuristic}
\end{align}}%
where the second equality follows from (\ref{equation:g-spec1}) and the third equality follows from (\ref{equation:m-fnq-flip}) and (\ref{equation:m-fnq-x}).  Note the change in the lower limit of the first summation symbol in the last equality.  

We consider the case $a=0$, $M\equiv 0 \pmod 3$, and recall that $T_{a,M}(q)$ will depend on an arbitrary $z$.  We write $j=t(M/3)+r$, where $0\le t\le 2$ and $0\le r\le (M/3)-1.$  The first summand in $(\ref{equation:DAMQ-heuristic})$ is then
{\allowdisplaybreaks \begin{align}
\frac{1}{M}&\sum_{j=0}^{M-1}m\big (q\zeta_M^{3j},q^3,-1\big )=\frac{1}{M}\sum_{t=0}^2\sum_{r=0}^{(M/3)-1}m\big (q\zeta_{M/3}^{t(M/3)+r},q^3,-1\big )\notag \\
&=\frac{3}{M}\sum_{r=0}^{(M/3)-1}m\big (q\zeta_{M/3}^{r},q^3,-1\big )\notag \\
&= m\big ( (-1)^{M+1} q^{M(M-1)/6}, q^{M^2/3}, z^{\prime} \big ) +\Psi_0^{M/3}(q,-1,z^{\prime};q^3),\label{equation:theta-last}
\end{align}}%
where the last equality follows from Theorem \ref{theorem:rootsof1} with $k=0$, $n=M/3$, $z=-1$.  The fourth summand is similar.  For the second summand,
{\allowdisplaybreaks \begin{align}
-\frac{1}{M}\sum_{j=0}^{M-1}\zeta_M^{-j}m\big (q\zeta_M^{-3j},q^3,-1\big )&=-\frac{1}{M}\sum_{t=0}^2\sum_{r=0}^{(M/3)-1}\zeta_M^{-(t(M/3)+r)}m\big (q\zeta_{M/3}^{-t(M/3)-r},q^3,-1\big )\notag \\
&=-\frac{1}{M}\Big ( \sum_{t=0}^2 \zeta_3^{-t}\Big )\sum_{r=0}^{(M/3)-1}\zeta_M^{-r}m\big (q\zeta_{M/3}^{-r},q^3,-1\big )= 0.
\end{align}}%
The third summand is similar and the result follows. So for $M\equiv 0 \pmod 3$ and $a=0$:
{\allowdisplaybreaks \begin{align}
D&(0,M,q)= 2m\big ( (-1)^{M+1} q^{M(M-1)/6}, q^{M^2/3}, z^{\prime} \big ) \notag\\
&\ \ \ \ \ +2\cdot \frac{ J_{M^2/3}^3} {\overline{J}_{0,3}}\cdot\sum_{t=0}^{(M/3)-1} \frac{q^{3\binom{t+1}{2}} j\big(-q^{3\binom{M/3+1}{2}+tM}/z^{\prime};q^{M^2/3}\big)
j\big ( (-1)^{M/3}q^{Mt+M/3} z^{\prime};q^{M^2/3}\big )}
{j\big( z^{\prime},(-1)^{M/3+1}q^{3\binom{M/3}{2}+M/3} z^{\prime}, (-1)^{M/3}q^{Mt+M/3};q^{M^2/3}\big)}\notag \\
&\ \ \ \ \ +\frac{1}{M}\sum_{j=1}^{M-1}\Big (1-\zeta_M^j\Big )  \frac{\zeta_M^{-j}J_1J_{3}^3j(-\zeta_M^{2j};q)}{j(\zeta_M^j;q)j(-q\zeta_M^{3j};q^3)j(-q^{2}\zeta_M^{3j};q^3)\overline{J}_{0,3}},\label{equation:Dfull-a-equal-0-and-3-div-M} 
\end{align}}%
where the theta quotients are from (\ref{equation:theta-last}) using (\ref{equation:msplit-theta}), and (\ref{equation:DAMQ-heuristic}) using (\ref{equation:theta-delta}).

We consider the case $a\equiv 1 \pmod 3$, $M\equiv 1 \pmod 3$, $M\neq 2a-1.$  Here $a=3b+1$, $M=3R+1$, so $a\equiv 3b-3R \pmod M$.  The first summand in $(\ref{equation:DAMQ-heuristic})$ is then
\begin{align}
\frac{1}{M}&\sum_{j=0}^{M-1}\zeta_M^{-aj} m\big (q\zeta_M^{3j},q^3,-1\big )
=\frac{1}{M}\sum_{j=0}^{M-1}\zeta_M^{3(\frac{M-1}{3}-\frac{a-1}{3})j}m\big (q\zeta_{M}^{3j},q^3,-1\big )\notag\\
&= (-1)^aq^{-a(a+1)/6-2M^2/3-2Ma/3-M/3}m\big ( (-1)^{M+1} q^{-M^2/2-Ma-M/2}, q^{3M^2}, -1 \big )\notag\\
&\ \ \ \ \  +\Psi_{(2M+a)/3}^M(q,-1,-1;q^3)\notag \\
&= (-1)^{M+a+1}q^{-a(a+1)/6-M^2/6+Ma/3+M/6} m\big ( (-1)^{M+1} q^{2M^2-M(3M-2a-1)/2}, q^{3M^2}, -1 \big )\notag \\
&\ \ \ \ \ +\Psi_{(2M+a)/3}^M(q,-1,-1;q^3), \label{equation:sumA}
\end{align}
where the second equality follows from Theorem \ref{theorem:rootsof1} with $k=(2M+a)/3$ and the third equality follows from (\ref{equation:m-fnq-flip}).  For the second summand in (\ref{equation:DAMQ-heuristic}) we write $a=3b+1$, $M=3R-2$, so $a+1\equiv 3b+2\equiv 3b+3R \pmod M$.  Here
\begin{align}
-\frac{1}{M}&\sum_{j=0}^{M-1}\zeta_M^{-(a+1)j} m\big (q\zeta_M^{-3j},q^3,-1\big )
=-\frac{1}{M}\sum_{j=0}^{M-1}\zeta_M^{-3(\frac{a-1}{3}+\frac{M+2}{3})j}m\big (q\zeta_{M}^{-3j},q^3,-1\big )\notag\\
&= (-1)^aq^{-a(a+1)/6-2M^2/3+2Ma/3+M/3}m\big ( (-1)^{M+1} q^{M^2-M(3M-2a-1)/2}, q^{3M^2}, -1 \big )\label{equation:sumB} \\
&\ \ \ \ \ -\Psi_{(2M-a-1)/3}^M(q,-1,-1;q^3), \notag
\end{align}
where the second equality follows from Theorem \ref{theorem:rootsof1} with $k=(2M-a-1)/3$.  
For the third summand in (\ref{equation:DAMQ-heuristic}) we write $a=3b+1$.  Here
{\allowdisplaybreaks \begin{align}
-\frac{1}{M}&\sum_{j=0}^{M-1}\zeta_M^{-(a-1)j} m\big (q\zeta_M^{3j},q^3,-1\big )
=-\frac{1}{M}\sum_{j=0}^{M-1}\zeta_M^{-3(\frac{a-1}{3})j}m\big (q\zeta_{M}^{3j},q^3,-1\big )\notag\\
&= (-1)^aq^{-a(a-1)/6}m\big ( (-1)^{M+1} q^{M^2+M(M-2a+1)/2}, q^{3M^2}, -1 \big )\notag \\
&\ \ \ \ \ -\Psi_{(a-1)/3}^M(q,-1,-1;q^3) \notag \\
&= (-1)^aq^{-a(a-1)/6}-(-1)^aq^{-a(a-1)/6}m\big ( (-1)^{M+1} q^{2M^2-M(M-2a+1)/2}, q^{3M^2}, -1\big )\label{equation:sumC}\\
&\ \ \ \ \ -\Psi_{(a-1)/3}^M(q,-1,-1;q^3), \notag
\end{align}}%
where the second equality follows from Theorem \ref{theorem:rootsof1} with $k=(a-1)/3$ and the third equality follows from (\ref{equation:m-fnq-x}) and then (\ref{equation:m-fnq-flip}).   For the fourth summand in (\ref{equation:DAMQ-heuristic}), we use $k=(M-a)/3$ to obtain
{\allowdisplaybreaks \begin{align}
\frac{1}{M}&\sum_{j=0}^{M-1}\zeta_M^{-aj} m\big (q\zeta_M^{-3j},q^3,-1\big )\notag \\
&= (-1)^{M+a}q^{-a(a-1)/6-M(M-2a+1)/6}m\big ( (-1)^{M+1} q^{M^2-M(M-2a+1)/2}, q^{3M^2}, -1 \big )\label{equation:sumD}\\
&\ \ \ \ \  +\Psi_{(M-a)/3}^M(q,-1,-1;q^3) \notag
\end{align}}%
With (\ref{equation:g-spec1}) in mind, combining (\ref{equation:sumA}) and (\ref{equation:sumB}) yields
{\allowdisplaybreaks \begin{align}
&(-1)^{a+1}q^{(2M^2-a(a+1))/6}g\big ( (-1)^{M+1} q^{M(3M-2a-1)/6},q^{M^2}\big ) \label{equation:finalA}\\
&\ \ \ \ \  -(-1)^{a+1}q^{(2M^2-a(a+1))/6}\Delta((-1)^{M+1}q^{M(3M-2a-1)/6};q^{M^2})\notag\\
&\ \ \ \ \ +\Psi_{(2M+a)/3}^M(q,-1,-1;q^3)-\Psi_{(2M-a-1)/3}^M(q,-1,-1;q^3). \notag
\end{align}}%
Combining (\ref{equation:sumC}) and (\ref{equation:sumD}) and using (\ref{equation:g-spec1}) yields
{\allowdisplaybreaks \begin{align}
 (-1)^{a}&q^{-a(a-1)/6}+(-1)^{M+a+1}q^{(M(M-2a+1)-a(a-1))/6}g\big ( (-1)^{M+1} q^{M(M-2a+1)/6},q^{M^2}\big )\label{equation:finalB}\\
 & -(-1)^{M+a+1}q^{(M(M-2a+1)-a(a-1))/6}\Delta((-1)^{M+1}q^{M(M-2a+1)/6};q^{M^2})\notag \\
 & -\Psi_{(a-1)/3}^M(q,-1,-1;q^3)  +\Psi_{(M-a)/3}^M(q,-1,-1;q^3). \notag
\end{align}}%
Combining (\ref{equation:finalA}) and (\ref{equation:finalB}) gives
\begin{align}
D(a,M,q)&=(-1)^{a}q^{-a(a-1)/6}\notag\\
&\ \ \ +(-1)^{M+a+1}q^{(M(M-2a+1)-a(a-1))/6}g\big ( (-1)^{M+1} q^{M(M-2a+1)/6},q^{M^2}\big )\notag \\
&\ \ \ +(-1)^{a+1}q^{(2M^2-a(a+1))/6}g\big ( (-1)^{M+1} q^{M(3M-2a-1)/6},q^{M^2}\big )+T_{a,M}(q),\notag
\end{align}
where
{\allowdisplaybreaks \begin{align}
T_{a,M}(q)&= -(-1)^{a+1}q^{(2M^2-a(a+1))/6}\Delta((-1)^{M+1}q^{M(3M-2a-1)/6};q^{M^2})\\
 & \ \ \ \ \ -(-1)^{M+a+1}q^{(M(M-2a+1)-a(a-1))/6}\Delta((-1)^{M+1}q^{M(M-2a+1)/6};q^{M^2})\notag \\
 &\ \ \ \ \ +\Psi_{(2M+a)/3}^M(q,-1,-1;q^3)-\Psi_{(2M-a-1)/3}^M(q,-1,-1;q^3) \notag\\
 & \ \ \ \ \ -\Psi_{(a-1)/3}^M(q,-1,-1;q^3)  +\Psi_{(M-a)/3}^M(q,-1,-1;q^3) \notag\\
 &\ \ \ \ \ +\frac{1}{M}\sum_{j=1}^{M-1}\zeta_M^{-aj}\Big (1-\zeta_M^j\Big ) \zeta_M^j\Delta(\zeta_M^j;q).\notag
\end{align} }%
Setting up for (\ref{equation:g-to-m}) instead of (\ref{equation:g-spec1}) in (\ref{equation:sumA})--(\ref{equation:sumD}) gives an alternate form:
\begin{align}
T_{a,M}(q)&=\Psi_{(2M+a)/3}^M(q,-1,q^{-M(3M-2a-1)/3};q^3)-\Psi_{(2M-a-1)/3}^M(q,-1,q^{M(3M-2a-1)/3};q^3) \notag\\
 & \ \ \ \ \ -\Psi_{(a-1)/3}^M(q,-1,q^{-M(M-2a+1)/3};q^3)  +\Psi_{(M-a)/3}^M(q,-1,q^{M(M-2a+1)/3};q^3) \notag\\
 &\ \ \ \ \ +\frac{1}{M}\sum_{j=1}^{M-1}\zeta_M^{-aj}\Big (1-\zeta_M^j\Big ) \zeta_M^j\Delta(\zeta_M^j;q).
\end{align}

%%%%%%%%%%%%%
%%%%%%%%%%%%%

\section{On Generalisations of Dyson's Rank-Differences}\label{section:BrO}

In \cite[Sections $7$, $5$]{Za}, Zagier obtains a statement slightly more general than \cite[Theorem $1.4$]{BrO2}.  He states that for all $M>0$ and all $a\in \mathbb{Z}/M\mathbb{Z}$ that the function (\ref{equation:Zag-def}) is a weight $1/2$ weakly holomorphic modular form for any arithmetic progression $\mathcal{A}\subset \mathbb{Z}$ not containing any number of the form $(1-h^2)/24$ with $h\equiv 2a\pm 1 \pmod{2M}$.  Zagier points out that this holds if $\mathcal{A}$ is the set of $n$  with $\Big ( \frac{1-24n}{\mathcal{P}}\Big )=-1$ for some prime $\mathcal{P}>3$, and in this case $\mathcal{A}$ is actually supported on a collection of arithmetic progressions.

It is straightforward to read such results from Theorem \ref{theorem:Theorem}.   We first focus on \cite[Theorem $1.4$]{BrO2}.  Take the case $M\equiv 0 \pmod 3$, $a\equiv 2 \pmod 3$ of Theorem \ref{theorem:Theorem}:
\begin{equation}
D(a,M,q)= (-1)^{a+1}\ q^{-a(a+1)/6}\ m((-1)^{M+1}\ q^{M(M-2a-1)/6},q^{M^2/3},z)+T_{a,M}(q),
\end{equation}
where $T_{a,M}(q)$ is a theta function which depends on the arbitrary $z$.  Using Theorem \ref{theorem:msplit-general-n} with $n=\mathcal{P}$ and $z'=-1$ (note that $n$ is from Theorem \ref{theorem:msplit-general-n} and not an $n\in\mathcal{A}$), we have
\begin{align}
D(a,M,q)&= (-1)^{a+1}\ q^{-a(a+1)/6} \sum_{r=0}^{\mathcal{P}-1} q^{-\frac{M^2}{3}\binom{r+1}{2}}\Big ( (-1)^{M}q^{\frac{M(M-2a-1)}{6}}\Big )^r \cdot  \notag \\
&\ \ \ \ \  \cdot \Big [ m(- q^{\frac{M^2}{3}\big (\binom{\mathcal{P}}{2}-\mathcal{P}r\big )}(  (-1)^{M}\ q^{\frac{M(M-2a-1)}{6}}  )^{\mathcal{P}},q^{\frac{(M\mathcal{P})^2}{3}},-1)\Big ]+T_{a,M}^{\prime}(q),
\end{align}
where $T_{a,M}^{\prime}(q)$ is $T_{a,M}(q)$ after having absorbed the theta functions which appear when using Theorem \ref{theorem:msplit-general-n}.  With the Appell--Lerch sum definition (\ref{equation:mxqz-def}) in mind, we can write
\begin{align}
D(a,M,q)&= (-1)^{Mr+a+1}\ q^{-a(a+1)/6} \sum_{r=0}^{\mathcal{P}-1} q^{-\frac{M^2}{3}\binom{r+1}{2}+\frac{rM(M-2a-1)}{6}}\sum_{m=m_r}^{\infty}a_r(m)q^{\mathcal{P}m} +T_{a,M}^{\prime}(q).\label{equation:mockness}
\end{align}
Multiplying by  $q^{-1/24}$ makes $T_{a,M}^{\prime}(q)$ a weight $1/2$ weakly holomorphic modular form.   We again remind the reader of the  various ways of producing new modular forms from old ones such as taking the twist by a Dirichlet character or restricting the index of summation of the Fourier expansion to terms that lie in certain arithmetic progressions.    So to ensure modularity, we only want the terms determined by the set $\mathcal{A}$ to be contributed from $T_{a,M}^{\prime}(q)$; in other words, we only have to worry about picking up `mockness'  from the $m(x,q,z)$ portion in (\ref{equation:mockness}) when we have for some $n\in \mathcal{A}$:  
\begin{equation}
-\frac{a(a+1)}{6}-\frac{M^2}{3}\binom{r+1}{2}+\frac{rM(M-2a-1)}{6}\equiv n \pmod {\mathcal {P}}.\label{equation:cond-1}
\end{equation}
Condition (\ref{equation:cond-1}) implies
\begin{equation}
1-24n\equiv (2rM+2a+1)^2 \pmod{\mathcal{P}}.
\end{equation}
The other subcases for the case $M\equiv 0 \pmod 3$ imply
\begin{equation}
1-24n\equiv (2rM+2a\pm 1)^2 \pmod{\mathcal{P}},
\end{equation}
where the sign depends on $a$.  So when 
\begin{equation}
\Big ( \frac{1-24n}{\mathcal{P}}\Big )=-1,
 \end{equation}
 then we have nothing to worry about.  The other cases $M\equiv 1,2\pmod 3$ are similar.  Here one needs Theorem  \ref{theorem:msplit-general-n} as well as (\ref{equation:g-full}).
%%%

%%%
We consider \cite[Theorem $1.1$]{BrOR}.  Here, $M\ge 5$ is prime.  For a given summand $q^mg(q^{Mk},q^{M^2})$, if $m\equiv d \pmod M$, it is easy to check that $1-24d\equiv (2a\pm 1)^2 \pmod {2M}$ where the sign depends on $M$ and $a.$  For example, in the first summand in the case $M\equiv 2 \pmod 3$, $a\equiv 1 \pmod 3$, the condition
\begin{equation}
\frac{M(M+2a-1)-a(a-1)}{6}\equiv d \pmod M
\end{equation}
implies
\begin{equation}
1-24d\equiv (2a-1)^2 \pmod {2M}.\label{equation:careful}
\end{equation}
If $\big ( \frac{1-24d}{M} \big ) =-1$,  we will never have an $m\equiv d \pmod M$, so any rank-difference can be taken.  If $\big ( \frac{1-24d}{M} \big ) =1$, then we may encounter $m$'s with $m\equiv d \pmod M$,  so one must be careful to avoid the rank values $a$ determined by (\ref{equation:careful}).

%%%%%%%%

\section{On $M_2$ rank-differences for partitions without repeated odd parts}\label{section:M2}

There are also variations of rank-differences.  Take for example the $M_2$ rank of Berkovitch and Garvan \cite{BG}, which is based on work of Dyson \cite{D2,D3}.   Here
\begin{align*}
M_2-\textup{rank}(\lambda)&:=\Big \lceil \tfrac{\ell (\lambda)}{2}\Big \rceil-\nu(\lambda)
\end{align*}
where $\ell (\lambda)$ is the largest part of $\lambda$ and $\nu(\lambda)$ is the number of parts of $\lambda$.  The $M_2$ rank of a partition is also the number of columns minus the number of rows in the $2$-modular diagram \cite{LO}, where, for example the partition $\lambda=(10,7,2,1)$ has $2$-modular diagram:
\[ \begin{array}{ccccc}
2&2&2&2&2\\
2&2&2&1\\
2\\
1
\end{array}
\]
Similar to the rank-differences of Dyson \cite{D1} and Atkin and Swinnerton-Dyer \cite{ASD}, Lovejoy and Osburn \cite[Theorems $1.1$, $1.2$]{LO} proved results on $M_2$ rank-differences for partitions without repeated odd parts.   

 Let $N_2(m,n)$ denote the number of partitions of $n$ without repeated odd parts whose $M_2$ rank is $m$.  We have the following nice generating function \cite[$(1.1)$]{LO}:
\begin{equation}
\sum_{\substack{n\ge0 \\ m\in \mathbb{Z}}}N_2(m,n)z^mq^n=\sum_{n=0}^{\infty}\frac{q^{n^2}(-q;q^2)_n}{(zq^2,q^2/z;q^2)_n}.
\end{equation}
Define $N_2(a,M;n)$ to be the number of partitions of $n$ without repeated odd parts whose $M_2$ rank is congruent to $a$ mod $M$, and define $p_o(n)$ to be the number of partitions of $n$ without repeated odd parts.   Furthermore we define the $M_2$ rank-difference
\begin{equation}
R_2(a,b,M,c,m):=\sum_{n= 0}^{\infty} \Big ( N_2(a,M;mn+c)-N_2(b,M;mn+c)\Big )q^n.\label{equation:R2-def}
\end{equation}

Lovejoy and Osburn showed \cite[Theorem $1.1$]{LO} for $M=m=3$, slightly rewritten:
{\begin{subequations}
\begin{equation}
R_2(0,1,3,0,3)=-1+3m(q^{15},q^{36},q^9)+3q^{-6}m(q^{-3},q^{36},q^{27})+{\overline{J}_{3,6}^2J_{1,4}}/{J_2^2}\label{equation:R_2-A}
\end{equation}
\begin{equation}
R_2(0,1,3,1,3)={\overline{J}_{3,6}J_{3,12}}/{J_2}\label{equation:R_2-B}
\end{equation}
\begin{equation}
R_2(0,1,3,2,3)={J_3J_{12}^2}/{J_4J_{1,6}}\label{equation:R_2-C}
\end{equation}
\end{subequations}}
We briefly sketch how we prove results for $M_2$ rank-differences similar to our Main Theorem.   We define the following:
\begin{equation}
D_2(a,M):=D_2(a,M,q)=\sum_{n=0}^{\infty}\Big ( N_2(a,M;n)-\frac{p_o(n)}{M}\Big )q^n.
\end{equation}
  From \cite{RLN, AM}, \cite[$(2.13)$]{Mo1}:
\begin{align}
\sum_{n= 0}^{\infty}\frac{(-1)^nq^{n^2}(q;q^2)_{n}}{(-x;q^2)_{n+1}(-q^2/x;q^2)_{n}}&=m(x,q,-1)+\frac{J_{1,2}^2}{2j(-x;q)},\label{equation:RLNid2}
\end{align}
so we can write
\begin{equation}
\sum_{\substack{n\ge0 \\ m\in \mathbb{Z}}}N_2(m,n)z^mq^n=\big (1-z\Big )\Big (m(-z,-q,-1) +\frac{\overline{J}_{1,2}}{2j(z,q)}\Big).
\end{equation}
It follows that 
\begin{equation}
D_2(a,M,q)=\sum_{j=1}^{M-1}\zeta_M^{-aj}\Big (1-\zeta_M^j\Big )\Big (m(-\zeta_M^j,-q,-1) +\frac{\overline{J}_{1,2}}{2j(\zeta_M^j,q)}\Big).\label{equation:D2-def}
\end{equation}

 Using (\ref{equation:D2-def}) and Theorem \ref{theorem:rootsof1} we obtain
\begin{theorem}\label{theorem:MT} If $0\le a<M$, then
\begin{align*}
D_2(a,M,q)= (-q)^{-\binom{a+1}{2}}m&\Big ( -(-q)^{\binom{M}{2}-Ma},(-q)^{M^2},z_0\Big)\\
&-(-q)^{-\binom{a}{2}}m\Big ( -(-q)^{\binom{M}{2}-M(a-1)},(-q)^{M^2},z_1\Big)+T_{a,M}(q),
\end{align*}
where $T_{a,M}(q)$ is a theta function which depends on $z_0$ and $z_1$.  Here $z_0$ and $z_1$ are arbitrary.
\end{theorem}
When multiplied by $q^{-1/8}$, the term $T_{a,M}(q)$ becomes a weight $1/2$ weakly holomorphic modular form.

Theorem \ref{theorem:MT} is useful for understanding results of Lovejoy and Osburn \cite[Theorems $1.1$, $1.2$]{LO}.   For the case $M=3$ we can use Theorem \ref{theorem:msplit-general-n} with $n=2$ to show
\begin{align*}
D_2(0,3)&=-1+2m(q^{15},q^{36},q^9)+2q^{-6}m(q^{-3},q^{36},q^{27})\\
&\ \ \ \ \ +\frac{2}{3}\cdot \frac{\overline{J}_{9,18}^2J_{3,12}}{J_6^2}+q\cdot \frac{2}{3}\cdot \frac{\overline{J}_{9,18}J_{9,36}}{J_{6}}+q^2\cdot \frac{2}{3}\cdot \frac{J_{9,36}^2}{J_{3,12}},\\
D_2(1,3)&=-m(q^{15},q^{36},q^9)-q^{-6}m(q^{-3},q^{36},q^{27})\\
&\ \ \ \ \ -\frac{1}{3}\cdot \frac{\overline{J}_{9,18}^2J_{3,12}}{J_6^2}-q\cdot \frac{1}{3}\cdot \frac{\overline{J}_{9,18}J_{9,36}}{J_{6}}-q^2\cdot \frac{1}{3}\cdot \frac{J_{9,36}^2}{J_{3,12}}.
\end{align*}
Identities (\ref{equation:R_2-A})--(\ref{equation:R_2-C}) are immediate.  For  $M=5$, we use Theorem \ref{theorem:msplit-general-n} with $n=2$ to show
 \begin{align*}
D_2&(0,5)=-1-2q^{-55}m(q^{-55},q^{100},q^{75})-2q^{-15}m(q^{-5},q^{100},q^{25})\\
&-\frac{6}{5}\cdot \frac{J_{5,50}^2J_{30,100}J_{40,100}}{J_{5}\overline{J}_{50,200}J_{100}}+q\cdot \frac{4}{5}\cdot \frac{J_{50}^4}{J_{10,50}J_{25}J_{100}}
 +q^2\cdot \frac{2}{5}\cdot \Big ( 2\cdot \frac{J_{25}^2J_{30,100}}{J_{10,25}J_{50}}-q^5\frac{J_{10,100}J_{25}^2}{J_{5,25}J_{50}}\Big )\\
 &-q^3\cdot \frac{2}{5}\cdot \frac{\overline{J}_{25,50}J_{25,100}}{J_{20,50}}
 +q^4\cdot \frac{2}{5}\cdot \frac{J_{10}J_{100}J_{15,50}^2}{J_{5}J_{30,100}J_{40,100}},\\
 D_2&(1,5)=-q^{-1}m(q^{35},q^{100},q^{25})-q^{-6}m(q^{15},q^{100},q^{25})\\
 &\ \ \ \ \ \ \ \ \ \  -1+q^{-55}m(q^{-55},q^{100},q^{75})+q^{-15}m(q^{-5},q^{100},q^{25})\\
&+\frac{4}{5}\cdot \frac{J_{5,50}^2J_{30,100}J_{40,100}}{J_{5}\overline{J}_{50,200}J_{100}}-q\cdot \frac{1}{5}\cdot \frac{J_{50}^4}{J_{10,50}J_{25}J_{100}}
 -q^2\cdot \frac{1}{5}\cdot \Big ( \frac{J_{25}^2J_{30,100}}{J_{10,25}J_{50}}-3\cdot q^5\frac{J_{10,100}J_{25}^2}{J_{5,25}J_{50}}\Big )\\
 &+q^3\cdot \frac{3}{5}\cdot \frac{\overline{J}_{25,50}J_{25,100}}{J_{20,50}}
 +q^4\cdot \frac{2}{5}\cdot \frac{J_{10}J_{100}J_{15,50}^2}{J_{5}J_{30,100}J_{40,100}},\\
D_2&(2,5)=q^{-1}m(q^{35},q^{100},q^{25})+q^{-6}m(q^{15},q^{100},q^{25})\\
&-\frac{1}{5}\cdot \frac{J_{5,50}^2J_{30,100}J_{40,100}}{J_{5}\overline{J}_{50,200}J_{100}}-q\cdot \frac{1}{5}\cdot \frac{J_{50}^4}{J_{10,50}J_{25}J_{100}}
 -q^2\cdot \frac{1}{5}\cdot \Big ( \frac{J_{25}^2J_{30,100}}{J_{10,25}J_{50}}+2\cdot q^5\frac{J_{10,100}J_{25}^2}{J_{5,25}J_{50}}\Big )\\
 &-q^3\cdot \frac{2}{5}\cdot \frac{\overline{J}_{25,50}J_{25,100}}{J_{20,50}}
 -q^4\cdot \frac{3}{5}\cdot \frac{J_{10}J_{100}J_{15,50}^2}{J_{5}J_{30,100}J_{40,100}}.
\end{align*}

It is straightforward to make generalisations on $M_2$ ranks similar to \cite[Theorem $1.4$]{BrO2} and \cite[Theorem $1.1$]{BrOR}

%%%%%%%%

%\enlargethispage{3\baselineskip}
\end{document}